\documentclass{article}

\usepackage{arxiv}

\usepackage[utf8]{inputenc} 
\usepackage[T1]{fontenc}    
\usepackage{hyperref}       
\usepackage{url}            
\usepackage{booktabs}       
\usepackage{amsfonts}       
\usepackage{nicefrac}       
\usepackage{microtype}      
\usepackage{lipsum}
\usepackage{graphicx}
\graphicspath{ {./images/} }
\usepackage{amssymb}
\usepackage{amsmath}

\title{Nuevo modelo para el dimensionamiento de lotes de pedidos en función del volumen de compra y deterioro temporal de los artículos}

\author{
 Margarita Miguelina Mieras \\
 Grupo de Investigación Sistemas Complejos (SiCo)\\
 Facultad Regional San Rafael \\
  Universidad Tecnol\'{o}gica Nacional\\
  Mendoza, CP 5600 \\
  \texttt{mmieras@frsr.utn.edu.ar} \\
   \And
 Tania Daiana Tobares \\
 Grupo de Investigación Sistemas Complejos (SiCo)\\
 Facultad Regional San Rafael \\
  Universidad Tecnol\'{o}gica Nacional\\
  Mendoza, CP 5600 \\
  \texttt{tanitobares@hotmail.com} \\
  \And
   Ricardo Raúl Palma\\
 Facultad de Ingenería \\
  Universidad Nacional de Cuyo\\
  Mendoza, CP 5500 \\
  \texttt{rpalma@uncu.edu.ar} \\
  \And
   Antonio José Ramirez-Pastor \\
 Instituto de Física Aplicada (INFAP)\\
 Universidad Nacional de San Luis\\
  CONICET\\
  San Luis, CP 5700 \\
  \texttt{antorami@gmail.com} \\
  \And
 Fabricio Orlando Sanchez Varretti \\
 Grupo de Investigación Sistemas Complejos (SiCo)\\
 Facultad Regional San Rafael \\
  Universidad Tecnol\'{o}gica Nacional\\
  Mendoza, CP 5600 \\
  \texttt{fsanchez@frsr.utn.edu.ar} \\
}

\begin{document}
\maketitle
\begin{abstract}
En esta investigación se presenta el desarrollo de un nuevo modelo de simulación para determinar el tamaño óptimo de los lotes de pedidos en la Planificación de Requerimientos de Materiales, en función del volumen de compra y deterioro temporal de los artículos. La novedad científica radica en el conteo exhaustivo de todas las estrategias de abastecimiento que se presentan a la hora de decidir cuándo y cuánta cantidad de materia prima y/o insumos adquirir gestionando simultáneamente múltiples factores. El modelo desarrollado permite obtener y visualizar la totalidad de soluciones factibles detectando un desafío significativo cuando se trabaja con horizontes de planificación de mayor tamaño. La metodología incluye el desarrollo de una ecuación matemática para calcular el costo total de todas las estrategias de abastecimiento, teniendo en cuenta descuentos por cantidad y plazo máximo permitido de los insumos en inventario. Los resultados muestran que el modelo analiza todo el espacio de búsqueda y halla la solución óptima. Se realiza la validación a través de la heurística búsqueda tabú, una técnica ampliamente reconocida en optimización. Se observa que si bien la heurística converge hacia el mínimo global requiere una carga computacional significativamente elevada. En contraste, el modelo desarrollado logra identificar el mínimo global con menor cantidad de cálculos, demostrando su eficiencia y precisión.
\end{abstract}


\section{Introducción}
El tamaño óptimo de los lotes de pedido en la Planificación de Requerimientos de Materiales (MRP, por sus siglas en inglés) es un concepto fundamental en la gestión de inventarios y la logística empresarial. La determinación de estos tamaños tiene implicaciones significativas en los costos operativos, la eficiencia del proceso de abastecimiento y la satisfacción del cliente. El tamaño del lote es uno de los problemas más importantes y también uno de los más difíciles en la planificación de la producción \cite{karimi}.

Existen enfoques clásicos que se encuentran en la literatura y abordan esta temática. Una de las técnicas fundamentales en la gestión de inventarios es el modelo de cantidad económica de pedido (EOQ, por sus siglas en inglés). Este modelo, desarrollado por Ford W. Harris en 1913 \cite{Harris} y posteriormente refinado por R.H. Wilson en 1934 \cite{Wilson}, establece que el tamaño óptimo de un lote de pedido se alcanza cuando el costo total de mantener inventario y el costo total de realizar pedidos son iguales. El Modelo de Balance Parcial del Periodo es una extensión del EOQ que permite la reposición de inventario en intervalos de tiempo discreto. En lugar de realizar pedidos de un tamaño fijo, se calcula la cantidad de pedido para cada periodo de tiempo basándose en el nivel de inventario al final del periodo anterior y la demanda esperada para el próximo periodo \cite{DeMatteis}. Otro de los modelos clásicos es el Algoritmo de Wagner y Whitin que es utilizado en la programación de producción y gestión de inventarios en entornos con capacidad limitada y demanda variable. El objetivo es minimizar los costos totales considerando los costos de producción, almacenamiento y mantenimiento de inventario, así como las restricciones de capacidad de producción \cite{Wagner}. Por su parte, el método heurístico Silver-Meal es utilizado para determinar el tamaño del lote de pedido en situaciones donde el comportamiento de la demanda es muy variable \cite{Silver}. Sin embargo, en situaciones en las que existen periodos de demanda cero no produce buenos resultados.

La temática sobre el problema de tamaño de lote sigue vigente, autores como Abolfazl, Seyed, Shekarabi y Karimi (2019) \cite{Abolfazl} se proponen optimizar el dimensionamiento de los lotes de pedido, mientras satisfacen restricciones estocásticas y calculan el número óptimo de lotes y el volumen de los mismos. Tobares, Mieras, Palma y Sanchez-Varretti (2023) \cite{Tobares} desarrollan un modelo de optimización para determinar el tamaño de lote óptimo considerando demanda constante y conocida.

Por otro lado, el deterioro de los insumos especialmente en productos perecederos añade una capa adicional de complejidad a la gestión de abastecimiento. Goyal y Giri (2001) \cite{Goyal2001} revisan modelos de inventario que consideran el deterioro y proponen métodos para optimizar los niveles en función de las tasas de deterioro. Chowdhury, Ghosh y Chaudhuri (2016) \cite{Chowdhury} investigan una política óptima de reabastecimiento para un artículo en deterioro permitiendo escasez de inventario. Mohr (2017) \cite{Mohr} propone encontrar decisiones de reabastecimiento óptimas sin tener toda la información de precios disponible. Barron (2018) \cite{Barron} presenta un problema de obsolescencia repentina de revisión continua estocástica similar a la EOQ. Huang, Jian y Tseng (2021) \cite{Huang}, mencionan que una posible solución sería analizar una cadena de suministro de dos niveles, compuesta por un proveedor y un minorista, donde los productos se deterioran con el tiempo y el minorista puede enfrentar situaciones de escasez. Asimismo, afirman que sería clave diseñar un mecanismo de coordinación basado en descuentos por cantidad que permita establecer de forma eficiente una política óptima de pedidos a largo plazo, con el objetivo de maximizar los beneficios de toda la cadena de suministro. Los descuentos por cantidad son incentivos ofrecidos por los proveedores para fomentar la compra en mayores volúmenes. Silver, Pyke y Peterson (1998) \cite{silver1998} exploraron la integración de descuentos por cantidad en el modelo EOQ. Este ajusta el tamaño del lote de pedido para aprovechar los descuentos, resultando en una reducción significativa de los costos totales cuando se alcanzan ciertos umbrales de pedido. Jhaveri y Gupta (2021) \cite{Jhaveri}demuestran que la coordinación entre vendedor y comprador, combinada con la implementación de descuentos por cantidad, genera una ganancia total adicional. Por esta razón, es esencial adoptar una estrategia de cadena de suministro coordinada que contemple estos descuentos.

Otros autores destacan la importancia de optimizar el inventario en sistemas con tasas de producción variables y características desafiantes de los productos, como el deterioro asociado a su vida útil limitada. Además, consideran que la demanda está influenciada por factores como el nivel de existencias, el estado de frescura del producto y su precio de venta. Para abordar este problema, emplean el método de Newton-Raphson como herramienta para obtener soluciones numéricas eficientes \cite{Tshinangi}. Según Rozi y Basri (2024) \cite{Rozi}, una coordinación eficiente entre varios proveedores, junto con la implementación de descuentos por cantidad, puede optimizar las ganancias dentro de la cadena de suministro. Además, desarrollan un modelo enfocado en integrar descuentos por cantidad en sistemas de inventario para productos con vidas útiles fijas, siguiendo la distribución de Weibul.

La literatura existente destaca la importancia de considerar múltiples factores en la determinación del tamaño óptimo de los lotes de pedido. Los descuentos por cantidad y el deterioro de los insumos son aspectos críticos que deben ser incorporados en los modelos de gestión de inventarios para reflejar de manera más precisa las condiciones reales del mercado. A pesar de los avances en la investigación sobre este tema, aún persisten desafíos importantes, especialmente en lo que respecta a la gestión simultánea de múltiples factores. Por lo tanto, en este trabajo se propone abordar esta brecha mediante el desarrollo de un modelo de optimización innovador que integre estos aspectos realizando un análisis exhaustivo de todas las estrategias de abastecimiento. Para validar este nuevo modelo se emplea la heurística de búsqueda tabú \cite{Glover}, \cite{GloverII} dada su amplia aceptación y su probada eficacia en problemas de optimización de investigación operativa \cite{Gopalakrishnan}, \cite{Li}. Al contrastar los resultados obtenidos con los de la heurística es posible evaluar la efectividad del modelo propuesto frente a una técnica consolidada, lo que brinda mayor contexto y solidez a los resultados. Además, su flexibilidad permite adaptar su implementación a las características específicas del problema, mientras que su capacidad para ajustar el número de iteraciones y las condiciones de parada facilita un análisis detallado de la eficiencia y calidad de las soluciones obtenidas. 

\section{Metodología}

En esta sección, se detalla el enfoque metodológico utilizado para abordar los objetivos de la investigación sobre el tamaño óptimo de lotes de pedido teniendo en cuenta factores como descuentos por cantidad y deterioro de los insumos.

\subsection{Totalidad de estrategias de abastecimiento }

El primer paso de la metodología de resolución es conocer todas las estrategias de abastecimiento que se pueden presentar al momento del aprovisionamiento de materiales. Para esto, se construyen matrices donde la cantidad de columnas esta dada por los periodos $i$ que constituyen el horizonte de planificación $N$ ($i=1,2,...,N$). El número de filas lo determina la totalidad de estrategias de abastecimiento $j$ donde $j$=1, 2,..,$2^{(N-1)}$. La Tabla 1 muestra una matriz genérica donde se observan estas características.

\begin{table}[h]
\begin{center}
\caption{Matriz A: representación de $j$ estrategias de abastecimiento para un horizonte de planificación de tamaño $N$.}
\begin{tabular}{|l|l|l|l|l|l|l|}
\hline
{Periodo $i$} & 1  & 2 & 3 & ... & N-1 & N \\ \hline
{Estrategia $j$} & $\alpha_1$ & $\alpha_2$ & $\alpha_3$ & ... & $\alpha_{(N-1)}$ & $\alpha_N$ \\ \hline
1 & 1 & 1 & 1 & ... & 1   & 1 \\ \hline
2 & 1  & 1  & 1  &  ... & 2  & 0  \\ \hline
$\vdots$ &   &   &   &     &     &   \\ \hline
j-1 & N-1  &  0 &  0 &   ...  &   0  &  1 \\ \hline
j & N & 0 & 0 & ... & 0   & 0 \\ \hline
\end{tabular}
\label{tabla1}
\end{center}
\end{table}

La cantidad de requerimientos de cada periodo se representa con $\alpha_i$ y el tamaño de lote de pedido esta dado por la coordenada $A[j,i]$.

\subsubsection{Estrategias de abastecimiento para un horizonte de planificación compuesto por 3 periodos}

Para ilustrar el desarrollo anterior se presenta un ejemplo sencillo. Para un horizonte de planificación $N=3$, se tiene un total de 4 estrategias de abastecimiento, ya que $j=2^{(N-1)}=2^{(3-1)}=4$. Se consideran 10, 20 y 15 unidades de requerimientos para los 3 periodos respectivamente, Tabla 2.

\begin{table}[h]
\begin{center}
\caption{Matriz A que representa $4$ opciones de pedido para un horizonte de planificación de tamaño $3$.}
\begin{tabular}{|l|lll|}
\hline
Periodo $i$     & \multicolumn{1}{l|}{1}  & \multicolumn{1}{l|}{2}  & 3  \\ \hline
Estrategia $j$ & \multicolumn{3}{c|}{$\alpha_{i}$}                    \\ \hline
                & \multicolumn{1}{l|}{10} & \multicolumn{1}{l|}{20} & 15 \\ \hline
1               & \multicolumn{1}{l|}{1}  & \multicolumn{1}{l|}{1}  & 1  \\ \hline
2               & \multicolumn{1}{l|}{1}  & \multicolumn{1}{l|}{2}  & 0  \\ \hline
3               & \multicolumn{1}{l|}{2}  & \multicolumn{1}{l|}{0}  & 1  \\ \hline
4               & \multicolumn{1}{l|}{3}  & \multicolumn{1}{l|}{0}  & 0  \\ \hline
\end{tabular}
\label{tabla2}
\end{center}
\end{table}

Como el tamaño del lote de pedido está dado por las coordenadas $A[j,i]$; la estrategia $j=1$ indica que se deben realizar tres pedidos en cada uno de los periodos, ya que: $A[1,1]$=1, $A[1,2]$=1, $A[1,3]$)=1. La cantidad solicitada dependera de los requerimientos $\alpha_i$, en este caso 10, 20 y 15, respectivamente.

Para $j=2$, los tamaños de los lotes son: $A[2,1]$=1, $A[2,2]$=2, $A[2,3]$=0. Esto indica que se deben realizar dos pedidos. El primero en el periodo 1 para cumplir con los requerientos de solo ese periodo, es decir, de 10 unidades. El segundo pedido de 35 unidades se lleva a cabo en el periodo 2 para cumplir con las necesidades de los periodos 2 y 3.

Para $j=3$, los tamaños de los lotes son: $A[3,1]$=2, $A[3,2]$=0, $A[3,3]$=1. Esto indica que se deben realizar dos pedidos. El primero en el periodo 1 para cumplir con los requerimientos de los periodos 1 y 2, es decir, de 30 unidades. El segundo pedido se lleva a cabo en el periodo 3 para cumplir con las necesidades de solo ese periodo, por lo tanto, se solicitan 15 unidades.

Para $j=4$, los tamaños de los lotes son: $A[4,1]$=3, $A[4,2]$)=0, $A[4,3]$=0. Esto indica que se realiza un solo pedido en el periodo 1 para cumplir con los requerimientos de todo el horizonte de planificación, es decir, el pedido es de 45 unidades.

Note que cuando la coordenada $A[j,i]$ es 0 no se realiza un pedido en el periodo $i$ analizado. El número de la coordenada $\neq$ de 0 indica la cantidad de periodos que son satisfechos con un solo pedido.

Observe que la cantidad de filas de estas matrices (estrategias de abastecimiento) responden a una fórmula exponencial lo que representa un desafio para su construcción cuando se trabaja con horizontes de planificación de mayor tamaño. La Figura 1 muestra este comportamiento y resalta algunos ejemplos: para un horizonte de planificación de 12 periodos el espacio de búsqueda estará compuesto por 2048 estrategias de pedido, para 24 periodos se deberán analizar 8388608 estrategias y para 30 periodos se estudiarán 5368870912 políticas de abastecimiento.

\begin{figure}[h]
\includegraphics[width=\textwidth]{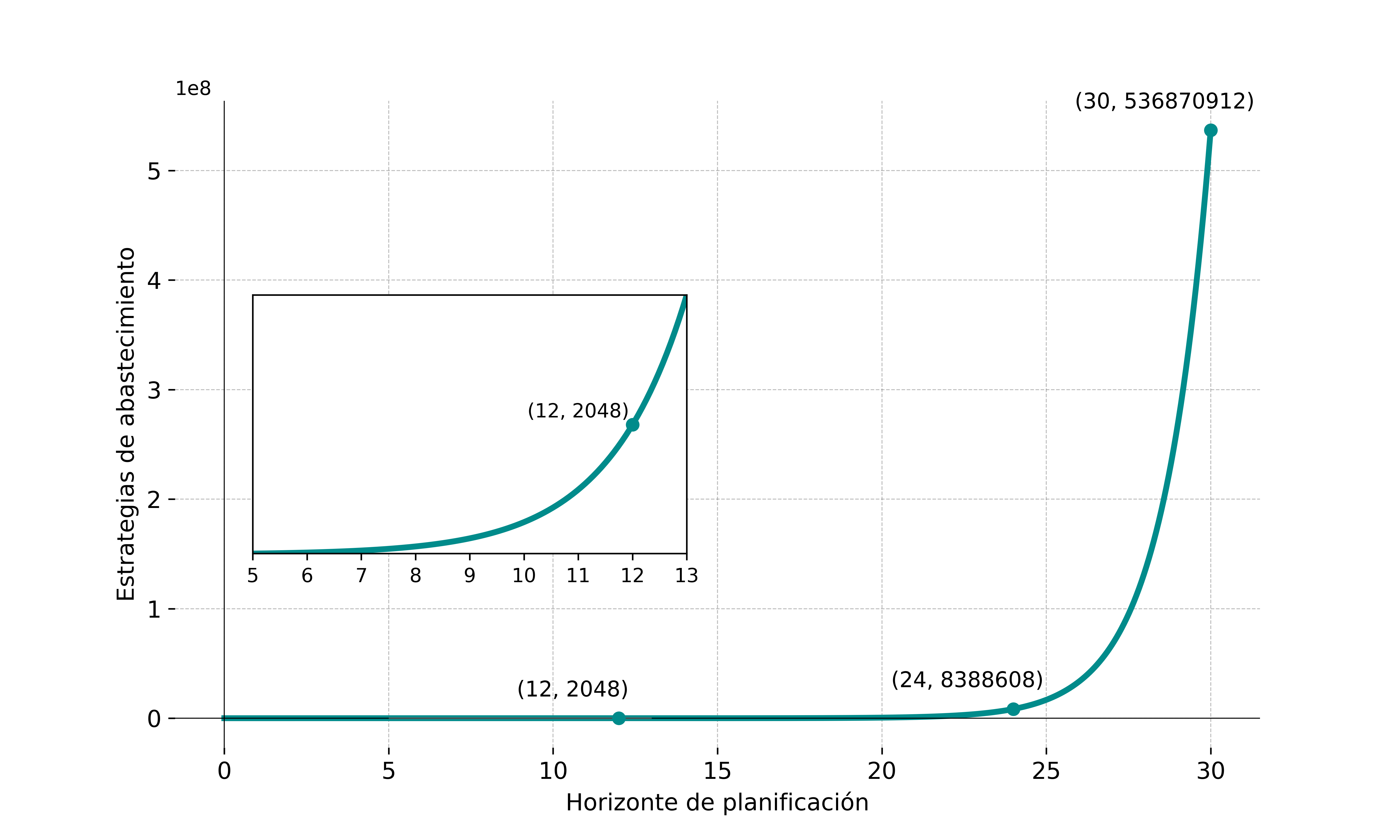}
\caption{Estrategias de abastecimiento vs. Horizonte de planificación.}
\label{fig.11}
\end{figure}

Ya desarrolladas las matrices que contienen todas las estrategias de abastecimiento, es decir, todo el espacio de búsqueda donde se halla la solución óptima se establece la función objetivo a minimizar.

\subsection{Modelo matemático}

El modelo matemático de este artículo se desarrolla sobre la base de las siguientes variables y supuestos.
\\

$Variables:$

$Cp$: Costo unitario de preparación de pedido.

$Cm$: Costo unitario de mantenimiento.

$Ca$: Costo unitario de adquisición.

$Cd$: Costo unitario de deterioro.

$\alpha$: Requerimientos.

$i$: Periodos.

N: Horizonte de planificación.

A[j,i]: Tamaño de lote.

$\mu$: plazo máximo permitido del insumo en inventario.

$q$: Cantidad de requerimiento a partir de la cual se aplica descuento.

$Np_{j}$: Cantidad de pedidos en la estrategia $j$
\\

$Supuestos:$

Se considera un solo elemento con una tasa de deterioro constante durante un horizonte de planificación conocido y finito de tamaño N.

La reposición se produce de manera instantánea.

No hay restricciones de espacio de almacenamiento.

No hay reparación o reemplazo de unidades deterioradas durante el horizonte de planificación. Los artículos serán retirados del almacén inmediatamente a medida que se deterioren.

No se permite escasez en los periodos $i$.
\\

De acuerdo con las variables y supuestos establecidos el modelo tiene cuatro costos generales, el costo total de pedido (CP), el costo total de mantenimiento (CM), el costo total de adquisión (CA) y el costo total de deterioro (CD), resultando el costo total (CT) de cada estrategia de abastecimiento:

\begin{equation}
CT = CP + CM + CA + CD
\end{equation}

El costo total de pedido (CP) resulta de la multiplicación del costo unitario de preparación de pedido por la cantidad de pedidos realizados:

\begin{equation}
CP = Cp·Np_{j}.
\end{equation}

El costo total de mantenimiento (CM) se obtiene multiplicando el costo unitario de mantenimiento por las cantidades almacenadas por el tiempo que estas permanece en inventario:

\begin{equation}
CM = Cm \left [\sum_{i=1}^{N}\sum_{i}^{l=(T[j,i]+i-1)} \alpha_l (l - i )\right ] \quad para \quad i \leqslant l.
\end{equation}

En la segunda sumatoria, cuando $i > l$, se tiene en cuenta el concepto de suma vacía que establece que en tal caso la sumatoria tiene el valor 0 \cite{Ingham}.\\

El costo total de adquisición (CA) se obtiene multiplicando el costo unitario de adquisición por las cantidades adquiridas:

\begin{equation}
CA = Ca \left [\sum_{i=1}^{N}\sum_{i}^{l=(T[j,i]+i-1)} \alpha_l\right ] \quad para \quad i \leqslant l.
\end{equation}

En la segunda sumatoria, cuando $i > l$, nuevamente se tiene en cuenta el concepto de suma vacía.

El costo de compra unitario $Ca$ está sujeto a descuentos por cantidad y responde a una función escalonada decreciente:


\begin{equation}
Ca(\alpha _{l})= Ca_{k} \quad si \quad q_{k-1}< \alpha _{l}\leq q_{k}, \quad Ca_{1}> Ca_{2}>...> Ca_{k}.
\end{equation}

Para un rango específico de cantidades ($\alpha _{l}$), el costo por unidad es un valor fijo ($Ca_{k}$).

$q_{k-1}$ y $q_{k}$ son los límites de ese rango. Si se compra una cantidad $\alpha _{l}$ que se encuentra entre estos dos límites, entonces pagas el costo $Ca_{k}$.

La relación $Ca_{1}> Ca_{2}> ...> Ca_{k}$ indica que los costos por unidad disminuyen a medida que se aumenta la cantidad comprada.\\

El costo total de deterioro (CD) se obtiene multiplicando el costo unitario de deterioro por las cantidades deterioradas por el tiempo que se excede el plazo máximo permitido del insumo en inventario:

\begin{equation}
CD = Cd \left [\sum_{i=1}^{N}\sum_{i}^{l=(T[j,i]+i-1)} \alpha_l ((l - i )-\mu)\right ] \quad para \quad i \leqslant l.
\end{equation}

Al igual que en los costos anteriores, en la segunda sumatoria, cuando $i > l$ se tiene en cuenta el concepto de suma vacía.

En este punto del trabajo se lograron desarrollar los cuatro términos generales que integran la fórmula objetivo de costo total, obteniendo:

\begin{multline}
    CT = Cp·Np_{j} + Cm \left [\sum_{i=1}^{N}\sum_{i}^{l=(T[j,i]+i-1)} \alpha_l (l - i )\right ] +\\
    Ca \left [\sum_{i=1}^{N}\sum_{i}^{l=(T[j,i]+i-1)} \alpha_l\right ] +  Cd \left [\sum_{i=1}^{N}\sum_{i}^{l=(T[j,i]+i-1)} \alpha_l ((l - i )-\mu)\right ] \quad para \quad i \leqslant l.   
\end{multline}

Con las bases conceptuales y metodológicas ya definidas, en la próxima sección se presentan los resultados obtenidos a partir de la aplicación del nuevo modelo de optimización propuesto. 

\section{Resultados}
En esta sección primero se detalla la verificación del modelo para asegurar que los cálculos y las relaciones matemáticas dentro de la fórmula son coherentes y no presentan errores lógicos o estructurales. Luego, se resuelve un ejemplo con un horizonte de planificación de 12 periodos simulando un escenario práctico. Este análisis permite explorar cómo el modelo gestiona las decisiones de abastecimiento en función de las variables definidas y evaluar su rendimiento en un contexto cercano a la realidad. Con ambas subsecciones se busca evidenciar la capacidad del modelo para analizar de forma exhaustiva todas las estrategias de abastecimiento y aportar soluciones óptimas para determinar el tamaño de los lotes de pedido.

\subsection{Verificación fórmula matemática}
En una primera instancia se realiza la verificación de la expresión costo total para comprobar que esta brinda resultados coherentes tanto desde el punto de vista numérico como dimensional.

En una primera instancia se realiza una verificación de la fórmula costo total para comprobar que esta brinda resultados coherentes tanto desde el punto de vista numérico como dimensional.

Retomando el ejemplo de la tabla 2 y considerando:

$Cp$=100 \$/pedido,   \quad $Cm$=1 \$/unidad.periodo, \quad $Ca_{1}$=5 \$/unidad, \quad $Ca_{2}$=4.50 \$/unidad, \quad $Cd$=10 \$/unidad.periodo, \quad $q$=30 unidades, \quad $\mu$=1 periodo, \quad $\alpha_{1}$=10 unidades, \quad $\alpha_{2}$=20 unidades, \quad $\alpha_{3}$=15 unidades, \quad $N$=3 periodos.

Se realizan los cálculos de los costos a través de un programa computacional de desarrollo propio en lenguaje $Python$ y se obtienen los resultados de cada estrategia de abastecimiento, Tabla 3.

\begin{table}[h]
\begin{center}
\caption{Costos totales}
\begin{tabular}{|l|l|l|l|l|l|}
\hline
 & CP & CM & CA & CD & CT \\ \hline
$j_{1}$ & 300 & 0 & 225 & 0 & 525 \\ \hline
$j_{2}$ & 200 & 15 & 207.5 & 0 & 422.5 \\ \hline
$j_{3}$ & 200 & 20 & 225 & 0 & 445 \\ \hline
$j_{4}$ & 100 & 50 & 202.5 & 150 & 502.5 \\ \hline
\end{tabular}
\label{tabla4}
\end{center}
\end{table}

Se corroboran los cálculos en pruebas de escritorio y son correctos, además se observa que exhiben coherencia lógica en sus resultados:

Los costos totales de preparación de pedido (CP) disminuyen a medida que se realiza menor cantidad de pedidos.

El costo total de mantenimiento (CM) es 0 en la estrategia de abastecimiento $j_1$ que implica la realización de pedidos en cada periodo sin opción a almacenamiento.

El costo total de adquisición (CA) varía en función del costo unitario de adquisición ($Ca$), el cual se ve condicionado por las cantidades adquiridas y si estas alcanzan el valor necesario $q$ para aplicar el descuento por cantidad.

El costo total de deterioro (CD) presenta un valor $\neq$ de 0 solo en la estrategia $j_4$ que es en la única oportunidad donde se supera el plazo máximo permitido del insumo en inventario, $\mu$.

Finalmente, el modelo brinda una solución óptima que establece que la estategía de abastecimiento con menor costo ($j_{2}$) implica realizar dos pedidos, el primero en el periodo $i_{1}$ de 10 unidades y el segundo y último pedido en el periodo $i_{2}$ de 35 unidades.

\subsection{Implementación del modelo en un horizonte de planificación de 12 periodos}
Considerando las siguientes variables:

$Cp$=100 \$/pedido,   \quad $Cm$=1 \$/unidad.pedido, \quad $Ca_{1}$=5 \$/unidad, \quad $Ca_{2}$=4.50 \$/unidad, \quad $Cd$=10 \$/unidad.periodo, \quad $q$=60 unidades, \quad $\mu$=4 periodos, \quad $\alpha_{1}$=10 unidades, \quad $\alpha_{2}$=20 unidades, \quad $\alpha_{3}$=10 unidades, \quad $\alpha_{4}$=35 unidades, \quad $\alpha_{5}$=40 unidades, \quad $\alpha_{6}$=40 unidades, \quad $\alpha_{7}$=35 unidades, \quad $\alpha_{8}$=10 unidades, \quad $\alpha_{9}$=5 unidades, \quad $\alpha_{10}$=5 unidades, \quad $\alpha_{11}$=35 unidades, \quad $\alpha_{12}$=40 unidades, \quad $N$=12 periodos.

Se realizan los cálculos de los costos a través de un programa computacional de desarrollo propio en lenguaje $Python$ y se obtienen los resultados de $2^{(N-1)}=2048$ estrategias de abastecimientos, Figura 2.

\begin{figure}[h]
\includegraphics[width=\textwidth]{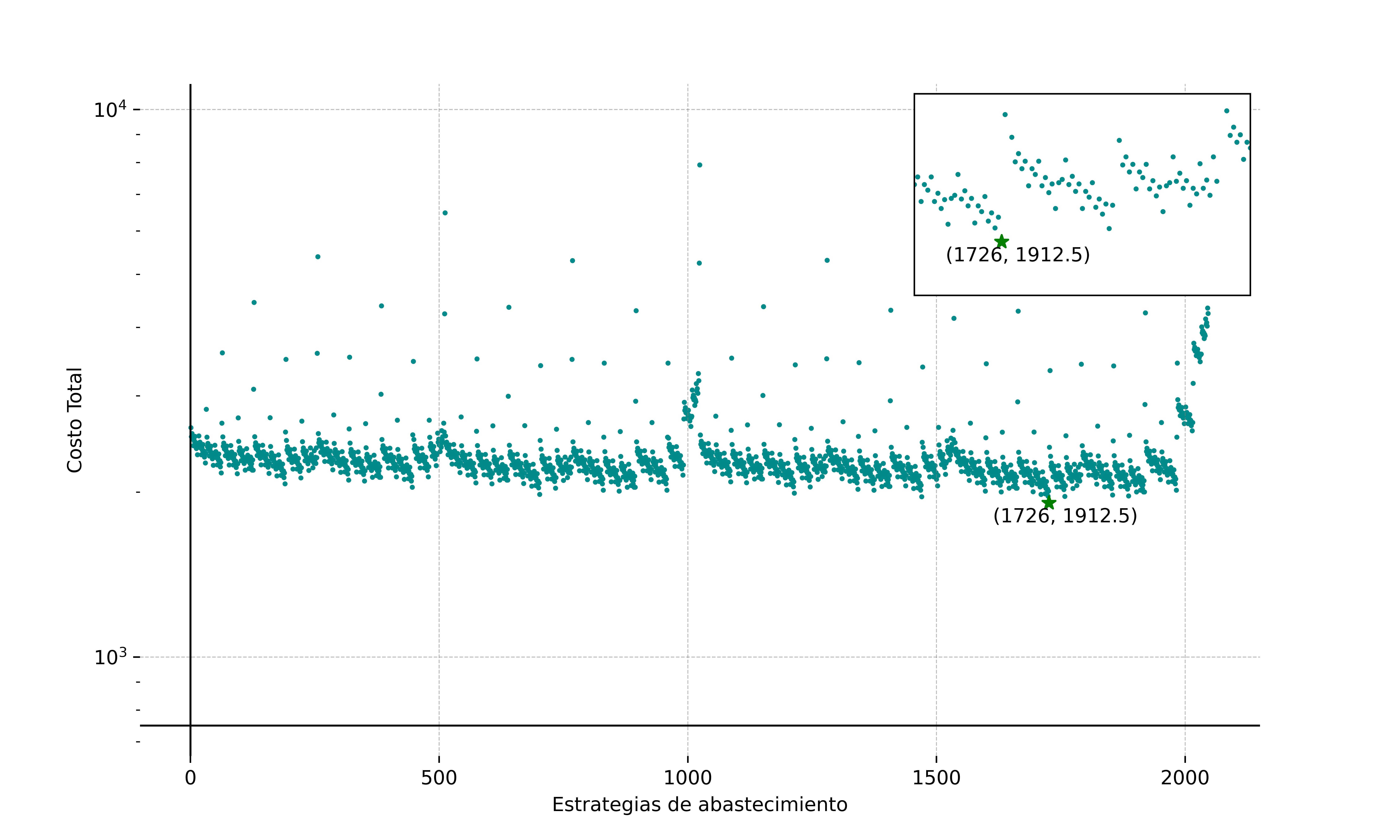}
\caption{Costo Total vs Estrategias de abastecimiento, N=12.} \label{fig.2}
\end{figure}

Se observa que la estrategia óptima es $j_{1726}$ ya que brinda el menor costo, $\$ 1912.5$. Al identificar esta estrategia en la matriz correspondiente a $N$=12, se obtiene $j_{1726}=[3, 0, 0, 2, 0, 5, 0, 0, 0, 0, 2, 0]$. Por lo tanto, se deben realizar cuatro pedidos, el primero en el periodo $i_{1}$ de 40 unidades, el segundo en el periodo $i_{4}$ de 75 unidades, el tercero en el periodo $i_{6}$ de 95 unidades y el último en el periodo $i_{11}$ de 75 unidades.

Los resultados muestran que el modelo brinda una visualización completa y detallada de todas las posibles estrategias de abastecimiento. Esto permite analizar el impacto de cada decisión en términos de costos totales, integrando factores clave como descuentos por cantidad y deterioro de insumos.

En esta sección se sientan las bases para comparar el desempeño del modelo desarrollado con alguna técnica reconocida, en este trabajo se ha optado por la heurística de búsqueda tabú.  Dicho método, cuya implementación se detalla en la sección Discusión, permitirá evaluar de manera precisa la efectividad y robustez del modelo propuesto frente a un referente reconocido en el campo de la optimización.

\section{Discusión}

Al utilizar la heurística búsqueda tabú para la validación se pretende demostrar la efectividad del modelo propuesto respecto a una técnica consolidada, dando mayor credibilidad y contexto a los resultados obtenidos. La flexibilidad de esta heurística permite implementar las especificidades del problema. Además, acepta el ajuste del número de iteraciones $k$ y condiciones de parada, lo que facilita el análisis de eficiencia.

Para dar comienzo al proceso de resolución se define la estructura de las soluciones, es decir de cada estrategia de pedido. Estas se representan mediante una serie binaria, donde:

0 implica no realizar pedido en ese periodo.

1 simboliza la realización de un pedido de las cantidades solicitadas en ese periodo y en los vecinos a la derecha con 0.

En el ejemplo presentado en la tabla 4 para un sistema de 12 periodos se realizan 5 pedidos:\\
1° pedido en $i_{1}$ solicitando $\alpha_{1}$.\\
2° pedido en $i_{2}$ solicitando $\alpha_{2}$, $\alpha_{3}$, $\alpha_{4}$, $\alpha_{5}$ y $\alpha_{6}$.\\
3° pedido en $i_{7}$ solicitando $\alpha_{7}$ y $\alpha_{8}$.\\
4° pedido en $i_{9}$ solicitando $\alpha_{9}$.\\
5° pedido en $i_{10}$ solicitando $\alpha_{11}$, $\alpha_{12}$ y $\alpha_{12}$.\\

\begin{table}[h]
\begin{center}
\caption{Estrategia de abastecimiento, N=12.}
\begin{tabular}{|l|l|l|l|l|l|l|l|l|l|l|l|}
\hline
1 & 2 & 3 & 4 & 5 & 6 & 7 & 8 & 9 & 10 & 11 & 12 \\ \hline
1 & 1 & 0 & 0 & 0 & 0 & 1 & 0 & 1 & 1 & 0 & 0 \\ \hline
\end{tabular}
\end{center}
\end{table}

La función objetivo a minimizar es la función de costo total presentada en el modelo propuesto. Esta heurística utiliza una lista tabú que ayuda a evitar la repetición de movimientos, lo que contribuye a escapar de mínimos locales. En este caso, se establece un periodo de permanencia en la lista tabú de $\tau$ = 5 iteraciones sucesivas. Esta configuración implica que, una vez añadido un movimiento a la lista, permanece allí durante cinco iteraciones evitando que el algoritmo explore soluciones recientemente visitadas y mejorando su capacidad de exploración en el espacio de soluciones.

El criterio de vecindad se establece mediante el intercambio aleatorio del valor de tres posiciones consecutivas en el vector de solución, intercambiando entre 0 y 1 en cada iteración. Este intercambio no afecta al primer elemento del vector, que permanece constante en 1, asegurando que se cumpla con las restricciones iniciales. Esta estrategia busca generar variabilidad en el vecindario, permitiendo al modelo explorar nuevas soluciones cercanas en cada iteración sin perder la consistencia en el valor inicial del vector.

Para que los resultados de la heurística sean comparables los valores de variables son los mismos que los utilizados en le punto 3.2.
\\ 

$Cp$=100 \$/pedido,   \quad $Cm$=1 \$/unidad.pedido, \quad $Ca_{1}$=5 \$/unidad, \quad $Ca_{2}$=4.50 \$/unidad, \quad $Cd$=10 \$/unidad.periodo, \quad $q$=60 unidades, \quad $\mu$=4 periodos, \quad $\alpha_{1}$=10 unidades, \quad $\alpha_{2}$=20 unidades, \quad $\alpha_{3}$=10 unidades, \quad $\alpha_{4}$=35 unidades, \quad $\alpha_{5}$=40 unidades, \quad $\alpha_{6}$=40 unidades, \quad $\alpha_{7}$=35 unidades, \quad $\alpha_{8}$=10 unidades, \quad $\alpha_{9}$=5 unidades, \quad $\alpha_{10}$=5 unidades, \quad $\alpha_{11}$=35 unidades, \quad $\alpha_{12}$=40 unidades, \quad $N$=12 periodos.\\

El proceso de cálculo se implementa mediante un programa computacional desarrollado en $Python$, diseñado específicamente para este estudio y adaptado a las especificaciones del problema. El algoritmo establece en forma aleatoria la estrategia inicial y realiza 10000 experimentos (z) para cada configuración de iteraciones establecida, con valores de $k$=10, 100, 1000, 5000 y 10000. En cada experimento guarda el costo mínimo que ha encontrado y finalmente calcula el costo promedio mínimo. Esta estructura permite una evaluación comparativa exhaustiva en distintas escalas de iteración, facilitando un análisis detallado del comportamiento y eficiencia en función de la cantidad de iteraciones. En la Figura 3 se muestra el costo promedio mínimo encontrado en cada grupo de experimentos y se comparan con el mínimo global encontrado a través del modelo presentado en este documento (línea roja continua).

\begin{figure}[h]
\includegraphics[width=\textwidth]{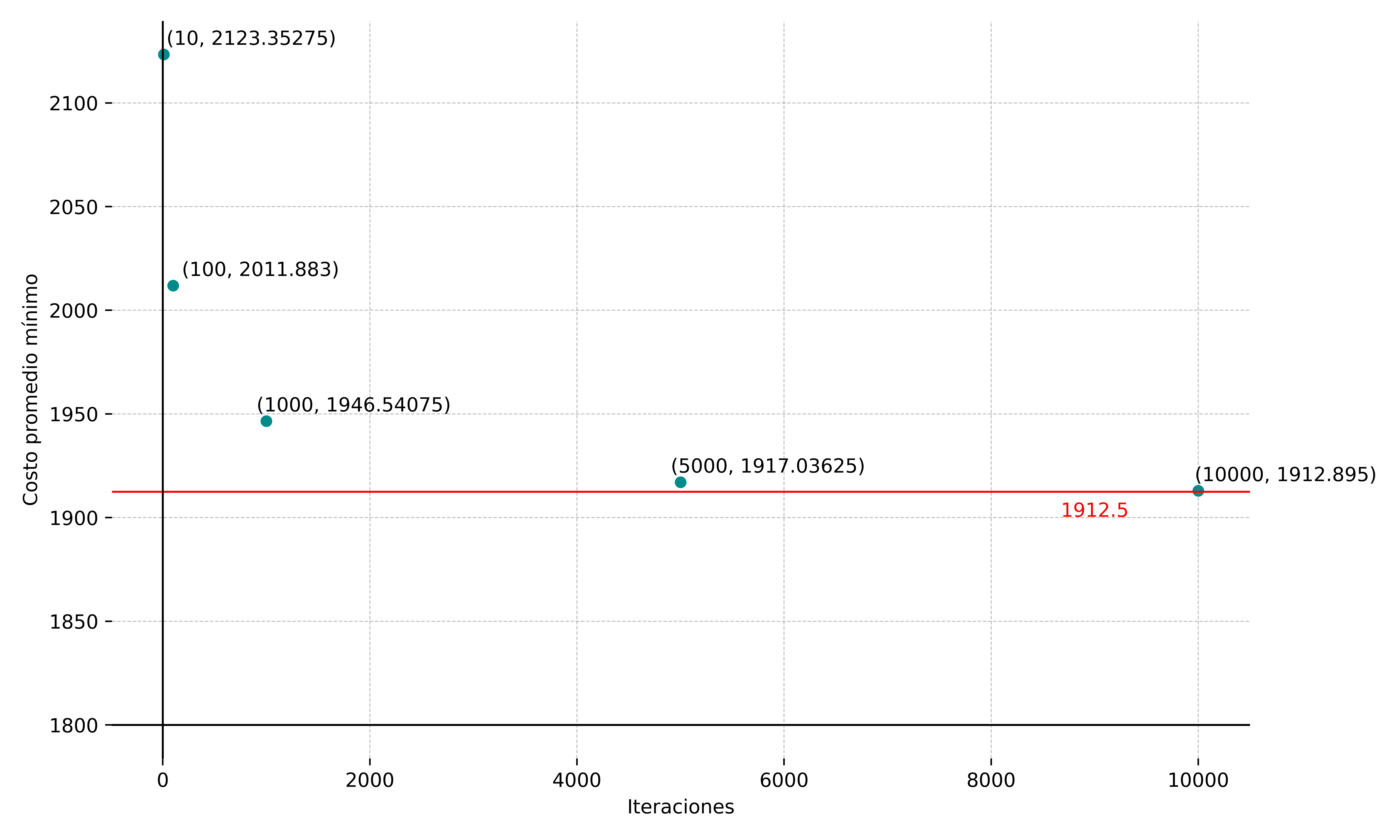}
\caption{Costo promedio mínimo (en 10000 experimentos) vs Iteraciones. Línea roja continua muestra mínimo global encontrado a través del nuevo modelo.} \label{fig.3}
\end{figure}

Se observa que a medida que se incrementa el número de iteraciones $k$ el costo promedio mínimo se acerca al mínimo absoluto encontrado mediante el modelo propuesto por los autores de este trabajo. Esto implica realizar cuatro pedidos, el primero en el periodo $i_{1}$ de 40 unidades, el segundo en el periodo $i_{4}$ de 75 unidades, el tercero en el periodo $i_{6}$ de 95 unidades y el último en el periodo $i_{11}$ de 75 unidades.

Sin embargo, no se puede asegurar que para un número dado de iteraciones la heurística encuentre el mínimo global. Si bien para k=10000 en la gran mayoría de los experimentos se encontró la solución óptima no sucedió el 100 \% de las veces. Para facilitar el análisis, el programa registra al final de cada experimento la serie binaria que representa la estrategia de menor costo encontrada al término de las k=10000 iteraciones. La Figura 4 muestra, de los $z$=10000 experimentos, cuantas veces la estrategia de menor costo emite un pedido en cada uno de los 12 periodos. Esto permite observar patrones en los periodos en los que es más frecuente que se realicen pedidos, ofreciendo así una visualización de las estrategias más eficientes en términos de costos acumulados.

\begin{figure}
\includegraphics[width=\textwidth]{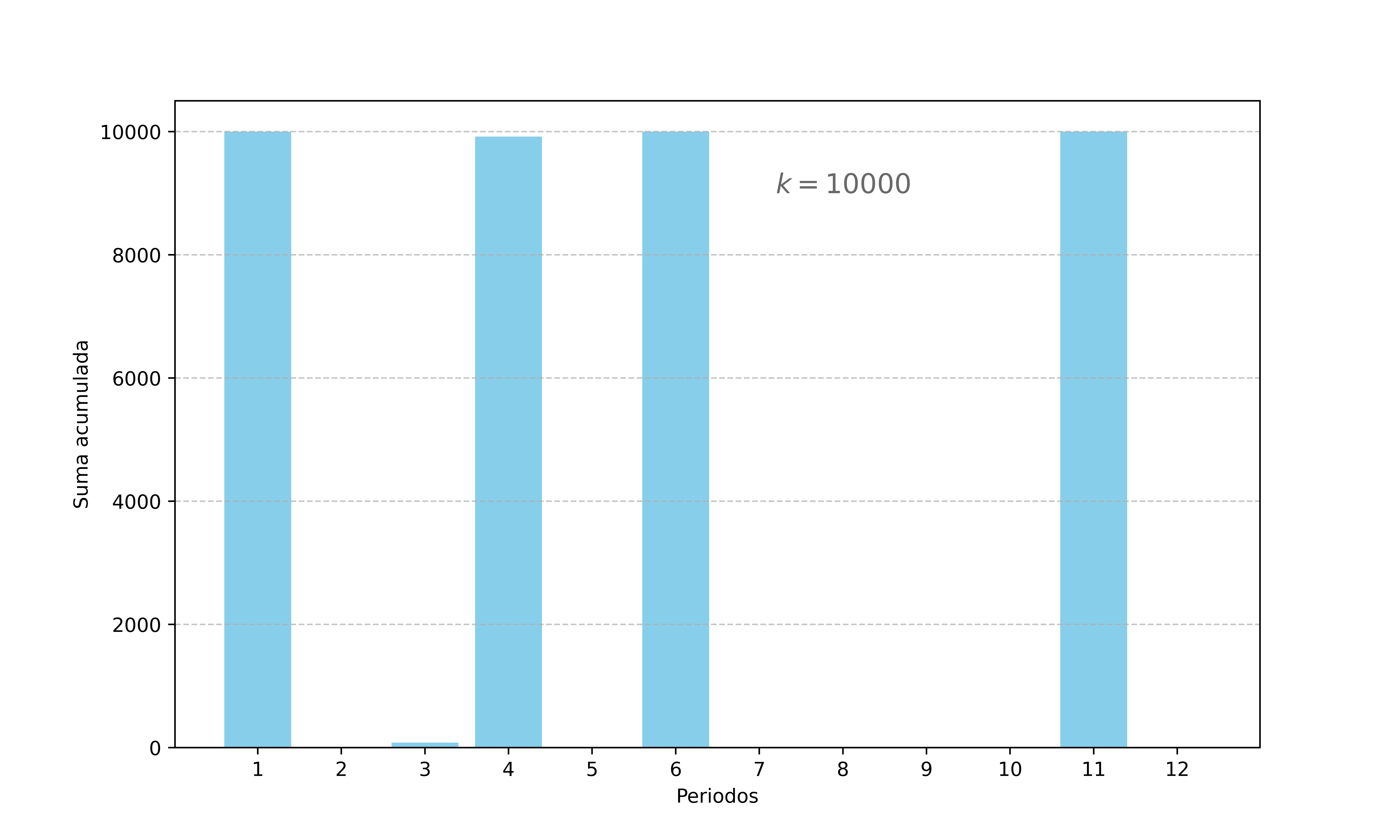}
\caption{Suma acumulada de los componentes de las estrategias de abastecimiento vs periodos} \label{fig.4}
\end{figure}

Este mismo análisis se puede realizar para las iteraciones $k$=10, 100, 1000 y 5000, Figura 5. En todas las soluciones de menor costo de cada experimento, el valor en el primer periodo es siempre alto, alcanzando el máximo posible de 10000. Esto indica que en el primer periodo siempre se realiza un pedido lo cual se debe a la condición inicial establecida en el modelo.

\begin{figure}
\includegraphics[width=\textwidth]{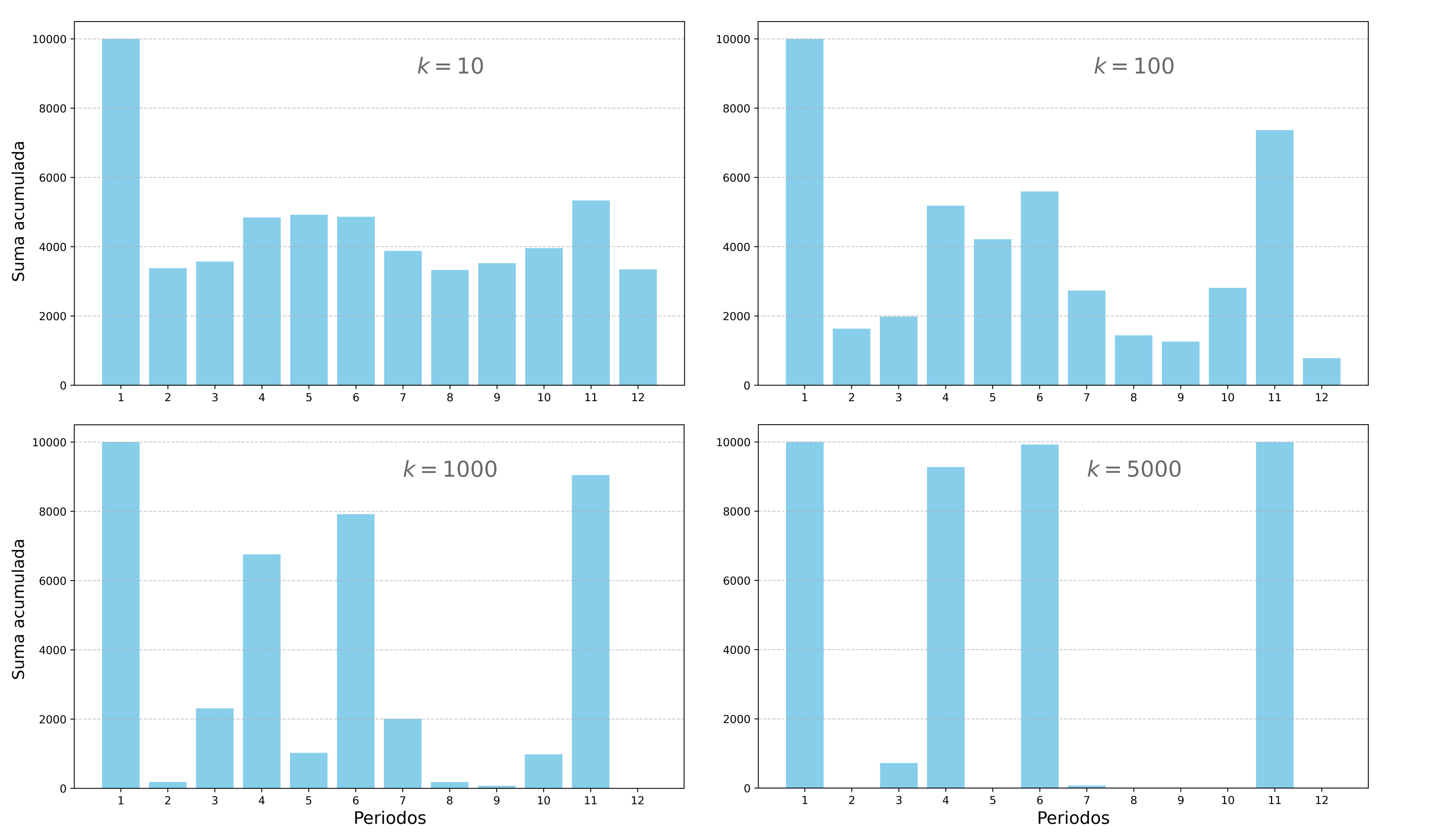}
\caption{Comparación suma acumulada de los componentes de las estrategias de abastecimiento vs periodos, para $k$=10, 100, 1000, 5000.} \label{fig.5}
\end{figure}

Para $k$=10 y $k$=100 se observa una distribución variada en los periodos restantes, lo que indica que a menor cantidad de iteraciones la heurística no muestra una tendencia hacia la solución de menor costo.

Para $k$=1000 y $k$=5000 los pedidos se centran en pocos periodos específicos (1, 3, 4, 6 y 11), y muchos periodos presentan valores cercanos a cero. Esto sugiere que a mayor cantidad de iteraciones el modelo converge hacia soluciones que favorecen ciertos periodos y omite otros indicando posibles patrones óptimos o soluciones repetitivas.

En el proceso de validación a través de la búsqueda tabú se observa que el nuevo modelo propuesto en este documento logra encontrar el mínimo global analizando sólo 2048 estrategias, lo que resalta su eficiencia en términos de cálculo y rapidez. En contraste, la heurística converge hacia ese mismo mínimo, pero requiere realizar millones de evaluaciones debido a la necesidad de explorar múltiples soluciones en el espacio de búsqueda. Esto evidencia que el modelo propuesto no solo es efectivo en encontrar la solución óptima, sino que también optimiza significativamente el uso de recursos computacionales.

\section{Conclusiones}

A partir del enfoque planteado a lo largo del documento, del análisis de los resultados obtenidos y la validación del modelo, se destacan las siguientes conclusiones:

i) En este estudio, se logra obtener y analizar exhaustivamente todas las soluciones factibles en el contexto del abastecimiento de materiales y la gestión de inventarios. Se desarrolla una fórmula matemática que permite calcular el costo total asociado a todas las estrategias de abastecimiento, considerando los supuestos y variables establecidas. Además, se comprueba que el nuevo modelo de optimización es capaz de explorar íntegramente el espacio de búsqueda y encontrar de manera eficiente la solución óptima.

ii) Un aspecto clave del modelo es la incorporación de factores como los descuentos por cantidad ofrecidos por los proveedores y el deterioro de los insumos, tomando en cuenta el tiempo máximo permitido para su almacenamiento. Estos elementos aportan una representación más precisa y realista de los costos y beneficios asociados a diferentes estrategias de abastecimiento, alineando el modelo con escenarios prácticos de la gestión en la cadena de suministro.

iii) Para validar la eficacia del modelo, se emplea la heurística de búsqueda tabú, una técnica consolidada en el campo de la optimización. Los resultados indican que, aunque la búsqueda tabú converge hacia el mínimo global, lo hace con un costo computacional significativamente mayor. En contraste, el modelo desarrollado alcanza la solución óptima con una cantidad considerablemente menor de evaluaciones lo que evidencia su alta eficiencia y precisión.

iv) Estos logros no solo resaltan los beneficios del modelo para gestionar el abastecimiento de materiales de manera eficiente, sino que también ofrecen una herramienta robusta y confiable para respaldar decisiones estratégicas. Tener acceso a la totalidad del espacio de soluciones permite identificar todas las rutas de acción posibles, ofreciendo un análisis exhaustivo que fundamenta sólidamente las decisiones en la gestión de la cadena de suministro.

v) Se detecta un desafío significativo en el cálculo del costo total de todas las estrategias de abastecimiento cuando se trabaja con horizontes de planificación más extensos, debido al crecimiento exponencial de las alternativas. En investigaciones futuras, se explorará la posibilidad de delimitar zonas de búsqueda con alta probabilidad de contener las mejores soluciones. Este enfoque tiene el potencial de reducir drásticamente los tiempos de cálculo y los costos computacionales, ampliando aún más la aplicabilidad del modelo en escenarios reales de gran escala.

Por último, se recomienda evaluar la integración del modelo con sistemas avanzados de apoyo a la decisión, como la inteligencia artificial o el aprendizaje automático, para fortalecer aún más su capacidad predictiva y adaptativa frente a entornos de alta incertidumbre y complejidad. Esto podría abrir nuevas oportunidades para su implementación en diversos sectores industriales, maximizando su impacto positivo en la optimización de recursos.

\end{document}